# MINIMAL SPANNING FORESTS

By Russell Lyons,[1] Yuval Peres[2] and Oded Schramm

*Indiana University, University of California, Berkeley and Microsoft Research*

Minimal spanning forests on infinite graphs are weak limits of minimal spanning trees from finite subgraphs. These limits can be taken with free or wired boundary conditions and are denoted FMSF (free minimal spanning forest) and WMSF (wired minimal spanning forest), respectively. The WMSF is also the union of the trees that arise from invasion percolation started at all vertices. We show that on any Cayley graph where critical percolation has no infinite clusters, all the component trees in the WMSF have one end a.s. In $\mathbb{Z}^d$ this was proved by Alexander [*Ann. Probab.* **23** (1995) 87–104], but a different method is needed for the nonamenable case. We also prove that the WMSF components are "thin" in a different sense, namely, on any graph, each component tree in the WMSF has $p_c = 1$ a.s., where $p_c$ denotes the critical probability for having an infinite cluster in Bernoulli percolation. On the other hand, the FMSF is shown to be "thick": on any connected graph, the union of the FMSF and independent Bernoulli percolation (with arbitrarily small parameter) is a.s. connected. In conjunction with a recent result of Gaboriau, this implies that in any Cayley graph, the expected degree of the FMSF is at least the expected degree of the FSF (the weak limit of uniform spanning trees). We also show that the number of infinite clusters for Bernoulli($p_u$) percolation is at most the number of components of the FMSF, where $p_u$ denotes the critical probability for having a unique infinite cluster. Finally, an example is given to show that the minimal spanning tree measure does not have negative associations.

**1. Introduction.** Suppose that i.i.d. random weights $\langle U(e); e \in \mathsf{E} \rangle$, uniform in $[0, 1]$, are assigned to the edges of a finite connected graph $G = (\mathsf{V}, \mathsf{E})$.

Received December 2004; revised July 2005.

[1]Supported in part by the Institute for Advanced Studies, Jerusalem, NSF Grants DMS-98-02663, DMS-01-03897 and DMS-02-31224, the Miller Institute, Berkeley.

[2]Supported in part by NSF Grants DMS-01-04073 and DMS-02-44479.

The minimal spanning tree determined by these weights is the spanning tree with minimum total weight; it can be obtained from $G$ by deleting every edge $e$ whose label $U(e)$ is maximal in some (simple) cycle. This construction has two analogues in an infinite graph $G$. The free minimal spanning forest FMSF in $G$ is obtained by deleting any edge with a label that is maximal in a cycle; the wired minimal spanning forest WMSF in $G$ is obtained by deleting any edge with a label that is maximal in an *extended* cycle, meaning a cycle or a bi-infinite simple path. The FMSF was studied by Alexander and Molchanov [6] in $\mathbb{Z}^2$ and by Alexander [5] in $\mathbb{Z}^d$. The WMSF is implicit in [5] (where it is shown that for $\mathbb{Z}^d$ it is the same as the FMSF); it was considered explicitly by Häggström [20], who studied the WMSF on trees, where it is generally different from the FMSF. Also, Aldous and Steele [3, 4] (Section 4.4 in the second paper) considered the wired minimal spanning forest (which they called simply the minimal spanning forest) in a weighted graph, in order to study the asymptotics of the analogous tree on the points of a Poisson process with weights given by Euclidean distances.

One reason to study these forests is their close connection to percolation; in fact, the WMSF is closely tied to critical bond percolation and invasion percolation, while the FMSF is related to percolation at the uniqueness threshold $p_u$ [as defined in (1.1)]. For example, the conjecture that nonamenable groups have an intermediate phase with infinitely many infinite clusters is equivalent to WMSF $\neq$ FMSF on such groups, as we shall see.

Another natural random forest in an infinite graph is the free uniform spanning forest (FSF), constructed by Pemantle [36] as a weak limit of uniform spanning trees from finite subgraphs. If the finite subgraphs are taken with wired boundary, then the wired uniform spanning forest (WSF) arises. These forests were studied in detail by Pemantle [36] and Benjamini, Lyons, Peres and Schramm (BLPS) [10] and have close connections to random walks and potential theory. A key theme of this paper will be to describe the striking analogies, and important differences, between the uniform and minimal spanning forests.

We now recall some terminology and state our main results. An *end* of an infinite tree $T$ is an equivalence class of infinite simple paths in $T$, where two paths are equivalent if they have a finite symmetric difference. Recall that a graph $G = (\mathsf{V}, \mathsf{E})$ is *transitive* if for every pair of vertices, there is an automorphism taking one to the other. A graph $G$ is *quasi-transitive* if the orbit space $\mathsf{V}/\mathrm{Aut}(G)$ is finite. A quasi-transitive graph $G$ is *unimodular* if $\mathrm{Aut}(G)$ is a unimodular group; see [9] for details. In particular, Cayley graphs are unimodular.

Recall that *Bernoulli*($p$) (*bond*) *percolation* on a graph $G$ is the random subgraph $G[p]$ that remains when each edge is independently kept with probability $p$ and deleted otherwise. Let $\theta(p) = \theta(p, G, v)$ be the probability that a fixed vertex $v$ of $G$ is in an infinite cluster of $G[p]$. As is customary, let $p_c = p_c(G)$ denote the infimum of those $p$ such that $G[p]$ has an infinite connected component a.s.; $p_c(G)$ is called the *critical probability* of $G$. Then $\theta(p) = 0$ for $p < p_c$ and $\theta(p) > 0$ for $p > p_c$ (assuming $G$ is connected). Bernoulli($p_c$) percolation is called *critical percolation*.

A very simple coupling among $\langle G[p]; p \in [0, 1] \rangle$, the WMSF and the FMSF is obtained by taking $G[p]$ to be the set of edges $e \in \mathsf{E}$ satisfying $U(e) < p$. This coupling, which we call the *standard coupling*, facilitates most of the connections between the WMSF, the FMSF and Bernoulli percolation.

We now describe our most interesting results.

THEOREM 1.1 (One end). *Let $G$ be a unimodular quasi-transitive graph. If $G$ is transitive, then the WMSF-expected degree of each vertex is* 2. *If $\theta(p_c, G) = 0$, then a.s. each component of the WMSF has one end.*

See Theorem 3.12 for the proof.

By [8], $\theta(p_c, G) = 0$ when $G$ is nonamenable (see below), quasi-transitive and unimodular. For $G = \mathbb{Z}^d$, Alexander [5] proved that there is one end per tree a.s. under the assumption that $\theta(p_c, \mathbb{Z}^d) = 0$. Kesten [27] and Hara and Slade [25] have shown that $\theta(p_c, \mathbb{Z}^d) = 0$ for $d = 2$ and $d \geq 19$, respectively, and this is widely believed to hold for every $d > 1$. Theorem 3.12 has an analogue for wired uniform spanning forests due to [36] and [10]; in the uniform case, no additional assumption like $\theta(p_c) = 0$ is needed. In any case, both theorems show that the wired spanning forests are "thin" in a certain sense.

Other aspects of "thinness" are the following. In [10], it was conjectured that the trees of the WSF on every graph are a.s. all recurrent for simple random walk. This was proved by Morris [33]. Our next theorem is an analogue of Morris' theorem for the WMSF.

THEOREM 1.2. *Let $G$ be an infinite graph. Then the WMSF $\mathfrak{F}_w$ satisfies $p_c(\mathfrak{F}_w) = 1$ a.s.*

See Theorem 3.20 for the proof of a slight improvement of this theorem.

On the other hand, the free spanning forests seem to be "thick" when not equal to the wired spanning forests. Thus, the next theorem shows that the FMSF is always "almost" connected:

THEOREM 1.3. *Let $G$ be any locally finite connected graph and $\varepsilon \in (0, 1)$. Let $\mathfrak{F}_f$ be a configuration of the FMSF and let $\omega$ be an independent copy of $G[\varepsilon]$. Then $\mathfrak{F}_f \cup \omega$ is connected a.s.*

See Theorem 3.22 for the proof. Note that for planar graphs, this follows from Theorem 1.2 by planar duality (see Theorem 4.1 below). In conjunction with recent work of Gaboriau [17], Theorem 1.3 allows us to compare the FSF and the FMSF:

COROLLARY 1.4. *For any transitive unimodular graph $G$, the FSF-expected degree of a vertex $v$ is at most the FMSF-expected degree of $v$.*

See Corollary 3.24 for the proof of an improvement of this theorem.

Obtaining a tight lower bound on the FSF-expected degree is a major open problem in the theory of uniform spanning forests, closely related to the possible equivalence of first $\ell^2$-Betti numbers and cost (see [16] and the proof of Corollary 3.24).

Let

$$(1.1) \quad p_{\mathrm{u}} = p_{\mathrm{u}}(G) := \inf\{p \in [0,1]; G[p] \text{ a.s. has a unique infinite cluster}\}.$$

PROPOSITION 1.5. *Let $G$ be an infinite graph. Under the standard coupling, the number of trees in the FMSF on $G$ is a.s. at least the number of infinite clusters of $G[p_{\mathrm{u}}]$.*

For a proof, see Proposition 3.13.

The preceding theorems show that in some ways, the WMSF is similar to critical Bernoulli percolation, while the FMSF is similar to Bernoulli percolation at the uniqueness threshold. Our next result holds for both minimal forests:

THEOREM 1.6. *Both measures WMSF and FMSF have a trivial tail $\sigma$-field on any graph.*

See Theorem 3.14 for the proof. The analogous result holds for the uniform spanning forests ([10] and [36]).

Naturally, we want to know when the free and wired minimal spanning forests are the same. We say that a graph $G$ has *almost everywhere uniqueness* (of the infinite cluster) if, for almost every $p \in (0,1)$ in the sense of Lebesgue measure, there is a.s. at most one infinite component in $G[p]$. The following extends Proposition 2.1 of [5]:

PROPOSITION 1.7. *On any connected graph $G$, we have* FMSF = WMSF *iff $G$ has almost everywhere uniqueness.*

See Proposition 3.6 for the proof.

For which $d$ is the FMSF in $\mathbb{Z}^d$ a tree? This is perhaps the most tantalizing open question about minimal spanning forests, and it has been answered only for $d = 2$; see [6]. A similar result of ours is as follows. Recall that an infinite graph is *nonamenable* if, for some $\delta > 0$ and all finite sets $F$ of vertices in $G$, the number of edges in the boundary of $F$ is at least $\delta|F|$.

PROPOSITION 1.8. *Let $G$ be a proper plane connected nonamenable graph with one end such that there is a group of homeomorphisms of the plane acting quasi-transitively on the vertices of $G$. Then the FMSF on $G$ is a.s. a tree.*

See Section 4 for the definition of a proper plane graph and Proposition 4.4 for the proof. This proposition applies to 1-skeletons (i.e., sets of vertices and edges) of tilings in the hyperbolic plane for which there is a group of hyperbolic isometries acting quasi-transitively on the tiles, and such that each tile is a finite and bounded hyperbolic polygon.

In the direction opposite to the proposition, it is not straightforward to produce any graph for which the FMSF is disconnected. An indirect method is to show that a certain transitive graph $G$ is *nontreeable*, that is, there is no automorphism-invariant measure on the space of spanning trees in $G$. It is known that Cayley graphs with Kazhdan's property $T$, as well as certain nonamenable products, are nontreeable (see [1] and [37]). As explained in the latter paper, in any nontreeable transitive graph, the number of components of the FMSF is almost surely infinite. Another indirect method is to apply Proposition 1.5 to transitive graphs where $G[p_\mathrm{u}]$ has infinitely many infinite clusters (such as the product of a regular tree of high degree and $\mathbb{Z}$; see [39]). This, however, does not yet provide graphs beyond the nontreeable graphs, since the only graphs for which this property of Bernoulli($p_\mathrm{u}$) percolation has been established are also known to be nontreeable. An example of a (nontransitive) graph where the FMSF has exactly two components a.s. is described in Example 6.1; it is not known whether this is possible in a transitive graph. Other examples and open questions on minimal spanning forests are presented in the final section. In Section 5 we give a simple method to calculate probabilities of spanning trees in finite graphs; we use this to give an example where the minimal spanning tree measure does not have negative associations, unlike the case of uniform spanning trees.

## 2. Background: minimal spanning trees on finite graphs.

We begin with a few definitions and some notation for graphs. A graph $G = (\mathsf{V}, \mathsf{E})$ is *locally finite* if the number of neighbors of each vertex is finite. We shall consider only such graphs. A *forest* is a graph with no cycles. A *tree* is a nonempty connected forest. A subgraph $H \subset G$ is *spanning* if $H$ contains all the vertices of $G$. A spanning tree or forest of $G = (\mathsf{V}, \mathsf{E})$ will usually be thought of as a subset of $\mathsf{E}$, since its vertex set is always $\mathsf{V}$. Given a graph $G = (\mathsf{V}, \mathsf{E})$, we let $2^\mathsf{E}$ denote the measurable space of all subsets of $\mathsf{E}$ with the Borel $\sigma$-field, that is, the $\sigma$-field generated by sets of the form $\{F \subset \mathsf{E}; e \in F\}$, where $e \in \mathsf{E}$.

Let $G = (\mathsf{V}, \mathsf{E})$ be a finite connected graph, and suppose that $U : \mathsf{E} \to \mathbb{R}$ is some injective function. The number $U(e)$ will then be referred to as the *label* of $e$. The labeling $U$ then induces a total ordering on $\mathsf{E}$, where $e < e'$ if $U(e) < U(e')$.

Define $T_U$ to be the subgraph whose vertex set is $\mathsf{V}$ and whose edge set consists of all edges $e \in \mathsf{E}$ whose endpoints cannot be joined by a path whose edges are strictly smaller than $e$. For the sake of completeness, we now prove that $T_U$ is a spanning tree. The largest edge in any cycle of $G$ is not in $T_U$ and, therefore, $T_U$ is a forest. If $\varnothing \neq A \subsetneq \mathsf{V}$, then the least edge of $G$ connecting $A$ with $\mathsf{V} \setminus A$ must belong to $T_U$, which shows that $T_U$ is connected. Thus, it is a spanning tree. In fact, among all spanning trees, $T_U$ has minimal edge label sum, $\sum_{e \in T} U(e)$.

DEFINITION 2.1 (*The minimal spanning tree*).   When $\langle U(e); e \in \mathsf{E} \rangle$ are independent uniform $[0, 1]$ random variables, the law of the corresponding spanning tree $T_U$ is called simply the *minimal spanning tree* (*measure*). It is a probability measure on $2^{\mathsf{E}}$.

There is an easy (and well-known) monotonicity principle for the minimal spanning tree measure, which is analogous to a similar principle for uniform spanning trees:

PROPOSITION 2.2 (Domination).   *Let $H$ be a connected subgraph of the finite connected graph $G$. Let $\mu_H$ and $\mu_G$ be the corresponding minimal spanning tree measures. Let $\mu_G^H(\mathcal{A}) := \mu_G[\{T; T \cap H \in \mathcal{A}\}]$ for every $\mathcal{A} \subseteq 2^{\mathsf{E}(H)}$. Then $\mu_H$ stochastically dominates $\mu_G^H$.*

A monotone coupling giving the stochastic domination is obtained by using the same labels on $\mathsf{E}(H)$, and independent labels on $\mathsf{E}(G) \setminus \mathsf{E}(H)$.

REMARK 2.3.   The following difference from the uniform spanning tree must be kept in mind. Given an edge, $e$, the minimal spanning tree measure on $G$ conditioned on the event not to contain $e$ need not be the same as the minimal spanning tree measure on $G \setminus e$, the graph $G$ with $e$ deleted; the simplest example is $G = K_4$, the complete graph on four vertices.

**3. Minimal spanning forests.**   We now present two natural extensions to infinite graphs of the notion of the minimal spanning tree. After the definitions, we give some partly novel deterministic facts about the forests, then proceed to our main probabilistic results.

Let $G = \langle \mathsf{V}, \mathsf{E} \rangle$ be an infinite connected graph and $U : \mathsf{E} \to \mathbb{R}$ be an injective labeling of the edges. Let $\mathfrak{F}_{\mathrm{f}} = \mathfrak{F}_{\mathrm{f}}(U) = \mathfrak{F}_{\mathrm{f}}(U, G)$ be the set of edges $e \in \mathsf{E}$ such that in every path in $G$ connecting the endpoints of $e$ there is at least one edge $e'$ with $U(e') \geq U(e)$. When $\langle U(e); e \in \mathsf{E} \rangle$ are independent uniform random variables in $[0, 1]$, the law of $\mathfrak{F}_{\mathrm{f}}$ (or sometimes, $\mathfrak{F}_{\mathrm{f}}$ itself ) is called the *free minimal spanning forest* on $G$ and is denoted by FMSF or FMSF$(G)$.

An *extended path* joining two vertices $a, b \in \mathsf{V}$ is either a simple path in $G$ joining them, or the union of a simple infinite path starting at $a$ and a disjoint simple infinite path starting at $b$. (The latter possibility may be considered as a simple path connecting $a$ and $b$ through $\infty$.) Let $\mathfrak{F}_{\mathrm{w}} = \mathfrak{F}_{\mathrm{w}}(U) = \mathfrak{F}_{\mathrm{w}}(U, G)$ be the set of edges $e \in \mathsf{E}$ such that in every extended path joining the endpoints of $e$ there is at least one vertex $e'$ with $U(e') \geq U(e)$. Equivalently, $\mathfrak{F}_{\mathrm{w}}(U)$ consists of those edges $e$ such that there is a finite set of vertices $W \subset \mathsf{V}$ where $e$ is the least edge joining $W$ to $\mathsf{V} \setminus W$. (If the endpoints of $e$ are $a$ and $b$, then $W$ is the vertex set of the component of $a$ or the component of $b$ in the set of edges smaller than $e$.) Again, when $U$ is chosen according to the product measure on $[0, 1]^{\mathsf{E}}$, we call $\mathfrak{F}_{\mathrm{w}}$ the *wired minimal spanning forest* on $G$. The law of $\mathfrak{F}_{\mathrm{w}}$ is denoted WMSF or WMSF(G).

Clearly, $\mathfrak{F}_{\mathrm{w}}(U) \subset \mathfrak{F}_{\mathrm{f}}(U)$. Note that $\mathfrak{F}_{\mathrm{w}}(U)$ and $\mathfrak{F}_{\mathrm{f}}(U)$ are indeed forests, since in every simple cycle of $G$, the edge $e$ with $U(e)$ maximal is present neither in $\mathfrak{F}_{\mathrm{f}}(U)$ nor in $\mathfrak{F}_{\mathrm{w}}(U)$. In addition, all the connected components in $\mathfrak{F}_{\mathrm{f}}(U)$ and in $\mathfrak{F}_{\mathrm{w}}(U)$ are infinite. Indeed, the least edge joining any finite vertex set to its complement belongs to both forests.

We shall now describe how $\mathfrak{F}_{\mathrm{f}}(U)$ and $\mathfrak{F}_{\mathrm{w}}(U)$ arise as limits of the minimal spanning tree on finite graphs. Consider an increasing sequence of finite, nonempty, connected subgraphs $G_n \subset G$, $n \in \mathbb{N}$, such that $\bigcup_n G_n = G$ and $\mathsf{E}(G_n) = (\mathsf{V}(G_n) \times \mathsf{V}(G_n)) \cap \mathsf{E}(G)$ for each $n$. For $n \in \mathbb{N}$, let $G_n^{\mathrm{w}}$ be the graph obtained from $G$ by identifying the vertices outside $G_n$ to a single vertex, then removing all resulting loops based at that vertex.

PROPOSITION 3.1. *Let $T_n(U)$ and $T_n^{\mathrm{w}}(U)$ denote the minimal spanning trees on $G_n$ and $G_n^{\mathrm{w}}$, respectively, that are induced by the labeling $U$. Then $\mathfrak{F}_{\mathrm{f}}(U) = \lim_{n \to \infty} T_n(U)$ and $\mathfrak{F}_{\mathrm{w}}(U) = \lim_{n \to \infty} T_n^{\mathrm{w}}(U)$. This means that for every $e \in \mathfrak{F}_{\mathrm{f}}(U)$, we have $e \in T_n(U)$ for every sufficiently large $n$, for every $e \notin \mathfrak{F}_{\mathrm{f}}(U)$, we have $e \notin T_n(U)$ for every sufficiently large $n$ and similarly for $\mathfrak{F}_{\mathrm{w}}(U)$.*

We leave the easy proof of this proposition to the reader.

It will be useful to make more explicit the comparisons that determine which edges belong to the two spanning forests. Define

$$Z_{\mathrm{f}}(e) = Z_{\mathrm{f}}^U(e) := \inf_{\mathcal{P}} \max\{U(e'); e' \in \mathcal{P}\},$$

where the infimum is over simple paths $\mathcal{P}$ in $G \setminus \{e\}$ that connect the endpoints of $e$; if there are none, the infimum is defined to be $\infty$. Thus, $\mathfrak{F}_{\mathrm{f}}(U) = \{e; U(e) \leq Z_{\mathrm{f}}(e)\}$. Similarly, define

$$Z_{\mathrm{w}}(e) = Z_{\mathrm{w}}^U(e) := \inf_{\mathcal{P}} \sup\{U(e'); e' \in \mathcal{P}\},$$

where the infimum is over extended paths $\mathcal{P}$ in $G \setminus \{e\}$ that join the endpoints of $e$. Again, if there are no such extended paths, then the infimum is defined to be $\infty$. Thus, $\{e; U(e) < Z_{\mathrm{w}}(e)\} \subseteq \mathfrak{F}_{\mathrm{w}}(U) \subseteq \{e; U(e) \leq Z_{\mathrm{w}}(e)\}$.

It turns out that there are also dual definitions for $Z_{\mathrm{f}}$ and $Z_{\mathrm{w}}$. In order to state these, recall that if $W \subseteq \mathsf{V}$, then the set of edges $\partial_{\mathsf{E}} W$ joining $W$ to $\mathsf{V} \setminus W$ is said to be a *cut*.

LEMMA 3.2. *For any injection $U : \mathsf{E} \to \mathbb{R}$ on any graph $G$, we have*

$$(3.1) \qquad Z_{\mathrm{f}}(e) = \sup_{\mathcal{C}} \inf \{U(e'); e' \in \mathcal{C} \setminus \{e\}\},$$

*where the supremum is over all cuts $\mathcal{C}$ that contain $e$. Similarly,*

$$(3.2) \qquad Z_{\mathrm{w}}(e) = \sup_{\mathcal{C}} \min \{U(e'); e' \in \mathcal{C} \setminus \{e\}\},$$

*where now the supremum is over all cuts $\mathcal{C}$ containing $e$ such that $\mathcal{C} = \partial_{\mathsf{E}} W$ for some finite $W \subset \mathsf{V}$.*

PROOF. We first verify (3.1). If $\mathcal{P}$ is a simple path in $G \setminus \{e\}$ that connects the endpoints of $e$, and $\mathcal{C}$ is a cut that contains $e$, then $\mathcal{C} \cap \mathcal{C} \neq \varnothing$, so $\max\{U(e'); e' \in \mathcal{P}\} \geq \inf\{U(e'); e' \in \mathcal{C} \setminus \{e\}\}$. This proves one inequality ($\geq$) in (3.1). To prove the reverse inequality, fix one endpoint $x$ of $e$, and let $W$ be the vertex set of the component of $x$ in $(G \setminus \{e\})[Z_{\mathrm{f}}(e)]$. Then $\mathcal{C} := \partial_{\mathsf{E}} W$ is a cut that contains $e$. Using $\mathcal{C}$ in the right-hand side of (3.1) yields the $\leq$ inequality in (3.1) and shows that the supremum there is achieved.

The $\geq$ inequality in (3.2) is proved in the same way as in (3.1). For the other direction, we dualize the above proof. Let $Z$ denote the right-hand side of (3.2), and let $W$ be the vertex set of the connected component of one of the endpoints of $e$ in the set of edges $e' \neq e$ such that $U(e') \leq Z$. We clearly have $U(e') > Z$ for each $e' \in \partial_{\mathsf{E}} W \setminus \{e\}$. Thus, by the definition of $Z$, the other endpoint of $e$ is in $W$ if $W$ is finite. The same argument applies with the roles of the endpoints of $e$ switched. Therefore, there is an extended path $\mathcal{C}$ connecting the endpoints of $e$ in $G \setminus \{e\}$ with $\sup\{U(e'); e' \in \mathcal{P}\} \leq Z$. This completes the proof of (3.2) and also shows that the infimum in the definition of $Z_{\mathrm{w}}(e)$ is attained. □

The *invasion tree* $T(v) = T_U(v)$ of a vertex $v$ is defined as the increasing union of the trees $\Gamma_n$, where $\Gamma_0 := \{v\}$ and $\Gamma_{n+1}$ is $\Gamma_n$ together with the least edge joining $\Gamma_n$ to a vertex not in $\Gamma_n$. (If $G$ is finite, we stop when $\Gamma_n$ contains $\mathsf{V}$.)

PROPOSITION 3.3. *Let $U : \mathsf{E} \to \mathbb{R}$ be an injective labeling of the edges of a locally finite graph $G = (\mathsf{V}, \mathsf{E})$. Then the union $\bigcup_{v \in \mathsf{V}} T_U(v)$ of all the invasion trees is equal to $\mathfrak{F}_{\mathrm{w}}(U)$.*

This is easily proved using the characterization of $\mathfrak{F}_w(U)$ as the set of all edges $e$ such that there is some finite $W \subset \mathsf{V}$ where $e$ is the minimal edge joining $W$ and $\mathsf{V} \setminus W$. (The details are left to the reader.)

The *invasion basin* $I(v)$ of a vertex $v$ is defined as the union of the subgraphs $G_n$, where $G_0 := \{v\}$ and $G_{n+1}$ is $G_n$ together with the lowest edge not in $G_n$ but incident to some vertex in $G_n$. Note that $I(v)$ has the same vertices as $T(v)$, but may have additional edges.

The following extends to general graphs a result proved in [6] (in $\mathbb{Z}^2$) and [5] (in $\mathbb{Z}^d$):

PROPOSITION 3.4. *Let $U : \mathsf{E} \to \mathbb{R}$ be an injective labeling of the edges of a locally finite graph $G = (\mathsf{V}, \mathsf{E})$. If $x$ and $y$ are vertices in the same component of $\mathfrak{F}_w(U)$, then the symmetric differences $I(x) \triangle I(y)$ and $T_U(x) \triangle T_U(y)$ are finite.*

PROOF. We give the proof only for $|I(x) \triangle I(y)| < \infty$, since the proof for $T_U(x) \triangle T_U(y)$ is essentially the same. It suffices to prove this when $e := [x, y] \in \mathfrak{F}_w(U)$. Consider the connected components $C(x)$ and $C(y)$ of $x$ and $y$ in $G[U(e)]$. Not both $C(x)$ and $C(y)$ can be infinite, since $e \in \mathfrak{F}_w(U)$. If both are finite, then invasion from each $x$ and $y$ will fill $C(x) \cup C(y) \cup \{e\}$ before invading elsewhere and, therefore, $I(x) = I(y)$ in this case. Finally, if, say, $C(x)$ is finite and $C(y)$ is infinite, then $I(x) = C(x) \cup \{e\} \cup I(y)$. $\square$

We begin our probabilistic results by recording the analogues of several results on uniform spanning forests from [10]:

PROPOSITION 3.5. *Let $G$ be a connected locally finite graph.*

(a) *If $G$ is amenable, then the average degree of vertices in both the free and wired minimal spanning forests on $G$ is a.s. 2.*

(b) *The free and wired minimal spanning forests on $G$ are the same if they have a.s. the same finite number of trees, or if the expected degree of every vertex is the same for both measures.*

(c) *The free and wired minimal spanning forests on $G$ are the same on any transitive amenable graph.*

(d) *If $\mathfrak{F}_w$ is connected a.s., or if each component of $\mathfrak{F}_f$ has a.s. one end, then $\mathrm{WMSF}(G) = \mathrm{FMSF}(G)$.*

(e) *If $G$ is unimodular and transitive with $\mathrm{WMSF}(G) \neq \mathrm{FMSF}(G)$, then a.s. the FMSF has a component with uncountably many ends, in fact, with $p_c < 1$.*

PROOF. The proofs are analogous to those of corresponding statements for uniform spanning forests in [10]. For (a), see Remark 6.1 there; for (b), see Remark 5.8 and Proposition 5.10 there; for (c), see Corollary 6.3 there; for (d), see Remark 5.9 there; and for (e), see Proposition 10.11 there.  □

Next, we characterize when the free and wired minimal spanning forests coincide.

PROPOSITION 3.6.  *On any connected graph $G$, we have* FMSF = WMSF *iff $G$ has almost everywhere uniqueness.*

PROOF.  Since $\mathfrak{F}_w \subset \mathfrak{F}_f$ and $\mathsf{E}$ is countable, FMSF $\neq$ WMSF is equivalent to the existence of an edge $e$ such that $\mathbf{P}[Z_w(e) < U(e) \leq Z_f(e)] > 0$. Let $A(e)$ be the event that the two endpoints of $e$ are in distinct infinite components of $(G \backslash e)[U(e)]$. Then $\{Z_w(e) < U(e) < Z_f(e)\} \subset A(e) \subset \{Z_w(e) \leq U(e) \leq Z_f(e)\}$. Consequently, $\mathbf{P}[A(e)] = \mathbf{P}[Z_w(e) < U(e) < Z_f(e)]$. Hence, FMSF $\neq$ WMSF is equivalent to the existence of an $e \in \mathsf{E}$ such that $\mathbf{P}[A(e)] > 0$. It is easy to see that almost everywhere uniqueness fails iff there is some $e \in \mathsf{E}$ with $\mathbf{P}[A(e)] > 0$.  □

COROLLARY 3.7.  *On any graph $G$, if almost everywhere uniqueness fails, then a.s.* WMSF *is not a tree.*

PROOF.  By Proposition 3.5(b), if WMSF is a tree a.s., then WMSF = FMSF.  □

We also obtain the following result of Häggström [22].

COROLLARY 3.8 (Equality on trees).  *If $G$ is a tree, then the free and wired minimal spanning forests are the same iff $p_c(G) = 1$.*

PROOF.  This is clear, since only at $p = 1$ can $G[p]$ have a unique infinite cluster a.s.; see [38] for this fact.  □

The issue of uniqueness in percolation is clarified by the following result. It was conjectured by Benjamini and Schramm [11] and proved by Häggström and Peres [23] under a unimodularity assumption, and by Schonmann [40] in general.

THEOREM 3.9 (Uniqueness monotonicity).  *Let $G$ be a locally finite connected graph with a quasi-transitive automorphism group. If there is a.s. a unique infinite cluster for $G[p]$, then the same holds for every $p' > p$.*

Thus, Proposition 3.6 shows that for quasi-transitive $G$, FMSF = WMSF iff $p_c = p_u$, which conjecturally holds iff $G$ is amenable. The argument of Burton and Keane [14] shows that, in fact, for a quasi-transitive amenable $G$ and every $p \in [0, 1]$, there is a.s. at most one infinite cluster in $G[p]$. This is slightly stronger than $p_c = p_u$, and provides another proof that for quasi-transitive amenable graphs, FMSF = WMSF [cf. Proposition 3.5(a), (b)].

In contrast to the WSF, the number of trees in the WMSF is not always an a.s. constant: see Example 6.2. On the other hand, the total number of ends of all trees in either forest is an a.s. constant, since it is a tail random variable and the tail $\sigma$-field is trivial, as we shall see in Theorem 3.14.

Although the number of trees in the WMSF can vary, we do know their essential supremum. Let $\alpha(x_1, \ldots, x_K)$ be the probability that $I(x_1), \ldots, I(x_K)$ are pairwise vertex-disjoint. The following theorem is analogous to Theorem 9.4 from [10]:

PROPOSITION 3.10. *Let $G$ be a connected graph. The WMSF-essential supremum of the number of trees is*

$$(3.3) \qquad \sup\{K; \exists x_1, \ldots, x_K \in \mathsf{V} \quad \alpha(x_1, \ldots, x_K) > 0\}.$$

This is obvious from the representation of WMSF as the union of invasion trees.

To analyze the number of ends in the trees of the WMSF, we shall use the Mass-Transport Principle, which was introduced into percolation theory by Häggström [21] and extended by BLPS [9]. See [9] for background. Also, the following lemma from [24] will be employed:

LEMMA 3.11. *Let $G$ be a locally finite quasi-transitive connected graph and $p > p_c(G)$. Then a.s. the invasion tree of each vertex of $G$ intersects an infinite cluster of $G[p]$.*

THEOREM 3.12 (One end). *Let $G$ be a unimodular quasi-transitive graph. If $G$ is transitive, then the WMSF-expected degree of each vertex is 2. If $\theta(p_c, G) = 0$, then a.s. each component of the WMSF has one end.*

PROOF. Fix a basepoint $o$. Let $e_1, e_2, \ldots$ be the edges in the invasion tree of $o$, in the order they are added. Suppose that $\theta(p_c) = 0$. Then $\sup_{n \geq k} U(e_n) > p_c$ for any $k$. By Lemma 3.11, $\limsup U(e_n) = p_c$. For each $k$ such that $U(e_k) = \sup_{n \geq k} U(e_n)$, the edge $e_k$ separates $o$ from $\infty$ in the invasion tree of $o$. It follows that the invasion tree of $o$ has a.s. one end. The same will be true for any finite connected union of invasion trees that contains $o$, since any such union agrees with the invasion tree of $o$, except for finitely many edges, by Proposition 3.4. Consequently, there is a well-defined special end for each component of $\mathfrak{F}_{\mathrm{w}}$ (viz., the end of any invasion tree contained in the component).

Suppose that $\mathrm{Aut}(G)$ acts transitively on the vertices of $G$. Orient each edge in $\mathfrak{F}_{\mathrm{w}}$ toward the special end of the component containing that edge. Then each vertex has precisely one outgoing edge. By the Mass-Transport Principle, it follows that the WMSF-expected degree of a vertex is 2. Since $\theta(p_{\mathrm{c}}) = 0$ when $G$ is non-amenable, this conclusion holds for all such $G$. It also holds for amenable $G$ by Proposition 3.5(a).

Combining the fact that the expected degree is 2 with Theorem 7.2 of [9], we deduce that a.s. each component of $\mathfrak{F}_{\mathrm{w}}$ has one or two ends. We say that a vertex $v$ is in the *future* of a vertex $u$ if $v$ can be reached from $u$ by following the oriented edges of $\mathfrak{F}_{\mathrm{w}}$.

Suppose that, with positive probability, the component of $o$ had two ends. Let the *trunk* of a component with two ends be the (unique) subgraph of it that would be isomorphic to $\mathbb{Z}$. Label the vertices of the trunk $x_n$ ($n \in \mathbb{Z}$), with $x_{n+1}$ in the future of $x_n$. Since $\theta(p_{\mathrm{c}}) = 0$, there would be an $\varepsilon > 0$ such that with positive probability $\sup_{n \in \mathbb{Z}} U([x_n, x_{n+1}]) > p_{\mathrm{c}} + \varepsilon$. By Lemma 3.11, $\limsup U([x_n, x_{n+1}]) = p_{\mathrm{c}}$. Thus, with positive probability there would be a largest $m \in \mathbb{Z}$ such that $U([x_m, x_{m+1}]) > p_{\mathrm{c}} + \varepsilon$. We could then transport mass 1 from each vertex in such a component to the vertex $x_m$, when this event occurs. The vertex $x_m$ would then receive infinite mass, contradicting the Mass-Transport Principle. This completes the proof in the case where $\mathrm{Aut}(G)$ acts transitively.

The reduction from the quasi-transitive setting to the transitive setting proceeds by mapping $\mathfrak{F}_{\mathrm{w}}$ to a forest over a set of vertices where there is a transitive group action. Let $\Gamma = \mathrm{Aut}(G)$, and denote by $\Gamma o := \{\gamma(o); \gamma \in \Gamma\}$ the orbit of $o$ under $\Gamma$. Let $k := \max_{v \in \mathsf{V}} \min_{o' \in \Gamma o} d(o', v)$, where $d(u, v)$ refers to the graph distance on $G$. Note that $k < \infty$, as $\Gamma$ acts quasi-transitively on $G$. For every vertex $v \in \mathsf{V}$, choose some vertex $o_v$ in $\Gamma o$ with $d(o_v, v) \leq k$. Among all possible choices of $o_v$, we choose uniformly at random, and make $\langle o_v; v \in \mathsf{V} \rangle$ independent and independent of the labeling $U$. Let $n$ be the number of vertices in $G$ within distance $k$ of $o$. For each $o' \in \Gamma o$, let $W_{o'} := \{v \in \mathsf{V}; o_v = o'\}$, and given $U$ and $\langle o_v; v \in \mathsf{V} \rangle$, choose an injection $\phi_{o'} : W_{o'} \to \{1, 2, \ldots, n\}$, uniformly at random among all possible injections, where the collection $\langle \phi_{o'}; o' \in \Gamma o \rangle$ is independent, given $U$ and $\langle o_v; v \in \mathsf{V} \rangle$. Let $\widetilde{\mathsf{V}} := \{1, 2, \ldots, n\} \times \Gamma o$, and let $\widetilde{\Gamma}$ denote the group of permutations of $\widetilde{\mathsf{V}}$ of the form $(j, o') \mapsto (\pi j, \gamma o')$, where $\pi$ is any permutation on $\{1, 2, \ldots, n\}$ and $\gamma \in \Gamma$. Then $\widetilde{\Gamma}$ is unimodular and acts transitively on $\widetilde{\mathsf{V}}$. Let $\widetilde{\mathfrak{F}}_{\mathrm{w}}$ be the image of $\mathfrak{F}_{\mathrm{w}}$ under the random injective map $\Phi : \mathsf{V} \to \widetilde{\mathsf{V}}$ given by $\Phi(v) := (\phi_{o_v}(v), o_v)$. Then $\widetilde{\mathfrak{F}}_{\mathrm{w}}$ is a forest on the vertex set $\widetilde{\mathsf{V}}$ and the law of $\widetilde{\mathfrak{F}}_{\mathrm{w}}$ is invariant under $\widetilde{\Gamma}$. The mass-transport argument used before now shows that the expected degree in $\widetilde{\mathfrak{F}}_{\mathsf{w}}$ of any vertex $\widetilde{v} \in \widetilde{\mathsf{V}}$ is $2\mathbf{P}[\widetilde{v} \in \Phi(\mathsf{V})]$. We may set $\widetilde{U}([\Phi(v), \Phi(u)]) := U([v, u])$ for every $[v, u] \in \mathsf{E}$. The above arguments then show that $\widetilde{\mathfrak{F}}_{\mathsf{w}}$ a.s. has precisely one end for every infinite connected component, which shows that the same is true for $\mathfrak{F}_{\mathsf{w}}$.  $\square$

The preceding result gives one relation between the WMSF and critical Bernoulli percolation. The next result gives a relation between the FMSF and Bernoulli($p_{\mathsf{u}}$) percolation.

According to Lemma 3.11, every component of WMSF($G$) intersects some infinite cluster of $G[p]$ for every $p > p_{\mathsf{c}}(G)$, provided $G$ is quasi-transitive. A comparable statement for FMSF($G$) holds for general graphs:

PROPOSITION 3.13. *Under the standard coupling, each component of* FMSF*($G$) intersects at most one infinite cluster of $G[p_{\mathsf{u}}]$. Thus, the number of trees in ($G$) is at least the number of infinite clusters in $G[p_{\mathsf{u}}]$. If $G$ is quasi-transitive with $p_{\mathsf{u}}(G) > p_{\mathsf{c}}(G)$, then each component of* FMSF*($G$) intersects exactly one infinite cluster of $G[p_{\mathsf{u}}]$.*

PROOF. Let $\langle p_j \rangle$ be a sequence satisfying $\lim_{j \to \infty} p_j = p_{\mathsf{u}}$ that is contained in the set of $p \in [p_{\mathsf{u}}, 1]$ such that there is a.s. a unique infinite cluster in $G[p]$. Let $\mathcal{P}$ be a finite simple path in $G$, and let $\mathcal{A}$ be the event that $\mathcal{P} \subset \mathfrak{F}_{\mathsf{f}}$ and the endpoints of $\mathcal{P}$ are in distinct infinite $p_{\mathsf{u}}$-clusters. Since a.s. for every $j = 1, 2, \ldots$ there is a unique infinite cluster in $G[p_j]$, a.s. on $\mathcal{A}$ there is a path joining the endpoints of $\mathcal{P}$ in $G[p_j]$. Because $\mathcal{P} \subset \mathfrak{F}_{\mathsf{f}}$ on $\mathcal{A}$, a.s. on $\mathcal{A}$ we have $\max_{\mathcal{P}} U \leq p_j$. Thus, $\max_{\mathcal{P}} U \leq p_{\mathsf{u}}$ a.s. on $\mathcal{A}$. On the other hand, $\max_{\mathcal{P}} U \geq p_{\mathsf{u}}$ a.s. on $\mathcal{A}$ since on $\mathcal{A}$, the endpoints of $\mathcal{P}$ are in distinct $p_{\mathsf{u}}$ components. This implies $\mathbf{P}[\mathcal{A}] \leq \mathbf{P}[\max_{\mathcal{P}} U = p_{\mathsf{u}}] = 0$, and the first statement follows.

The second sentence follows from the fact that every vertex belongs to some component of $\mathfrak{F}_{\mathsf{f}}$. Finally, the third sentence follows from Lemma 3.11 and the fact that invasion trees are contained in the wired minimal spanning forest, which, in turn, is contained in the free minimal spanning forest.  $\square$

We now prove a result that shows precisely shared behavior for both minimal spanning forests:

THEOREM 3.14. *Both measures WMSF and FMSF have a trivial tail $\sigma$-field on any graph.*

For our proof, we need the following strengthening of Theorem 5.1(i) of [5]:

LEMMA 3.15. *Let $G$ be any infinite locally finite graph with distinct fixed labels $U(e)$ on its edges. Let $\mathfrak{F}$ be the corresponding free or wired minimal spanning forest. If the label $U(e)$ is changed at a single edge $e$, then the forest changes at most at $e$ and at one other edge [an edge $f$ with $U(f) = Z_f(e)$ or $Z_w(e)$, resp.]. More generally, if $\mathfrak{F}'$ is the forest when labels only in $K$ are changed, then $|(\mathfrak{F} \triangle \mathfrak{F}') \setminus K| \leq |K|$.*

PROOF. Consider first the free minimal spanning forest. Suppose the two values of $U(e)$ are $u_1$ and $u_2$, with $u_1 < u_2$. Let $\mathfrak{F}_1$ and $\mathfrak{F}_2$ be the corresponding free minimal spanning forests. Then $\mathfrak{F}_1 \setminus \mathfrak{F}_2 \subseteq \{e\}$. Suppose that $f \in \mathfrak{F}_2 \setminus \mathfrak{F}_1$. Then there must be a path $\mathcal{P} \subset G \setminus \{e\}$ joining the endpoints of $e$ and containing $f$, such that $U(f) = \max_{\mathcal{P}} U > u_1$. Suppose that there were a path $\mathcal{P}' \subset G \setminus \{e\}$ joining the endpoints of $e$, such that $\max_{\mathcal{P}'} U < U(f)$. Then $\mathcal{P} \cup \mathcal{P}'$ would contain a cycle containing $f$ but not $e$, on which $f$ has the maximum label. This would contradict $f \in \mathfrak{F}_2$. Therefore, $Z_f(e) = U(f)$. Since the labels are distinct, there is at most one such $f$.

For the WMSF, the proof is the same, only with "extended path" replacing "path" and "$Z_w(e)$" replacing "$Z_f(e)$."

The second conclusion in the lemma follows by induction from the first. □

PROOF OF THEOREM 3.14. Let $\mathcal{F}(K)$ be the $\sigma$-field generated by $U(e)$ for $e \in K$. Let $A$ be a tail event of $2^{\mathsf{E}}$. Let $\phi : [0,1]^{\mathsf{E}} \to 2^{\mathsf{E}}$ be the map that assigns the (free or wired) minimal spanning forest to a configuration of labels. (Actually, $\phi$ is defined only on the configurations of distinct labels.) We claim that $\phi^{-1}(A)$ lies in the tail $\sigma$-field $\bigcap_{K \text{ finite}} \mathcal{F}(\mathsf{E} \setminus K)$. This implies the desired result by Kolmogorov's 0–1 law. Indeed, for any finite set $K$ of edges and any two labelings $\omega_1, \omega_2$ that differ only on $K$, we know by Lemma 3.15 that $\phi(\omega_1)$ and $\phi(\omega_2)$ differ at most on $2|K|$ edges, whence either both $\omega_i$ are in $\phi^{-1}(A)$ or neither are. In other words, $\phi^{-1}(A) \in \mathcal{F}(\mathsf{E} \setminus K)$. □

Assume, as usual, that $\langle U(e); e \in \mathsf{E} \rangle$ are uniform i.i.d. in $[0, 1]$. If $G$ is a transitive graph and $I(v)$ is the invasion basin (not just the tree) of a vertex $v$, then $p_c(I(v)) = 1$ a.s. (To prove this, let $p > p_c(G)$. Lemma 3.11 implies that $I(v) \cap G[p]$ is finite a.s. However, Bernoulli($p'$) percolation on $G[p]$ has the same law as $G[p'p]$. Thus, a.s. $p_c(I(v)) \geq p_c(G[p]) = p_c(G)/p$.)

In fact, a stronger result is true in greater generality. Define the *invasion basin of infinity* to be the set of edges $[x, y]$ such that there do not exist disjoint infinite simple paths from $x$ and $y$ consisting only of edges $e$ satisfying $U(e) < U([x, y])$, and denote the invasion basin of infinity by $I(\infty) = I^U(\infty)$. Note that

$$I(\infty) \supset \bigcup_{v \in \mathsf{V}} I(v) \supset \mathfrak{F}_w(U).$$

For an edge $e$, define

$$Z_\infty^U(e) := Z_\infty(e) := \inf_{\mathcal{P}} \sup\{U(f); f \in \mathcal{P} \setminus \{e\}\},$$

where the infimum is over bi-infinite simple paths that contain $e$; if there is no such path $\mathcal{P}$, define $Z_\infty(e) := 1$. Because $U(e)$ and $Z_\infty(e)$ are independent and $U(e)$ is a continuous random variable, we have $U(e) \neq Z_\infty(e)$ a.s. for every $e$. Thus, a.s. $e \in I(\infty)$ iff $U(e) < Z_\infty(e)$.

THEOREM 3.16. *Let $G = (\mathsf{V}, \mathsf{E})$ be a bounded degree graph. Then $p_c(I(\infty)) = 1$ a.s. Therefore, $p_c(\mathfrak{F}_w) = 1$ a.s.*

To prove this, we begin with the following lemma that will provide a coupling between percolation and invasion that is different from the usual one we have been working with:

LEMMA 3.17. *Let $G = (\mathsf{V}, \mathsf{E})$ be a locally finite infinite graph and $\langle U(e); e \in \mathsf{E} \rangle$ be i.i.d. uniform $[0,1]$ random variables. Conditioned on $I(\infty)$, the random variables*

$$\frac{U(e)}{Z_\infty(e)} \qquad (e \in I(\infty))$$

*are i.i.d. uniform $[0,1]$ random variables.*

REMARK 3.18. At first sight, this lemma may seem obvious; however, the proof requires some care, as the parameters $Z_\infty(e)$ are defined in terms of the edge variables $U(e)$. The proof given below circumvents this dependence. It relies on the following elementary observation:

Let $\langle U_i; 1 \le i \le k \rangle$ be a random vector distributed uniformly in $[0,1]^k$, and let $\langle Z_i; 1 \le i \le k \rangle$ be an independent random vector with an arbitrary distribution in $(0,1]^k$. Then given $U_i < Z_i$ for all $1 \le i \le k$, the conditional law of the vector $\langle U_i/Z_i; 1 \le i \le k \rangle$ is uniform in $[0,1]^k$.

We leave the justification to the reader.

PROOF OF LEMMA 3.17. Let $A \subset \mathsf{E}$ be a finite set. It suffices to prove that, conditional on $I^U(\infty)$, on the event $A \subset I^U(\infty)$, the random variables $\langle U(e)/Z_\infty(e); e \in A \rangle$ are i.i.d. uniform in $[0,1]$. Define $\widetilde{U}(e) := 0$ for $e \in A$ and $\widetilde{U}(e) := U(e)$ for $e \notin A$, and let $Z_A^U := Z_\infty^U \upharpoonright A$ denote the restriction of $Z_\infty^U$ to $A$.

*Claim:* The symmetric difference between the event $[U \upharpoonright A < Z_A^{\widetilde{U}}]$ and the event $[A \subset I^U(\infty)]$ has probability 0.

The first event is contained in the second event because $Z_A^{\widetilde{U}}(e) \le Z_A^U(e)$ for all $e \in A$. For the converse, assume that $A \subset I^U(\infty)$. Consider any bi-infinite simple path $\mathcal{P}$. If $e \in A \cap \mathcal{P}$, then $U(e) < Z_\infty^U(e) \le \sup\{U(e'); e \ne e' \in \mathcal{P}\}$. Hence, for every such $\mathcal{P}$,

$$\sup_{\mathcal{P}} U = \sup_{\mathcal{P} \backslash A} U = \sup_{\mathcal{P} \backslash A} \widetilde{U} = \sup_{\mathcal{P}} \widetilde{U}.$$

Therefore, $Z_A^U = Z_A^{\widetilde{U}}$ and $I^U(\infty) = I^{\widetilde{U}}(\infty)$, provided that $Z_\infty^U(e) \ne U(e)$ for all $e \in \mathsf{E}$, which holds a.s. Hence, $A \subset I^U(\infty)$ implies $U \restriction A < Z_A^{\widetilde{U}}$ a.s., which verifies the claim.

By the claim (and the observation in Remark 3.18), conditioned on $\widetilde{U}$ and $A \subset I^U(\infty)$, the random variables $\langle U(e)/Z_A^{\widetilde{U}}(e); e \in A \rangle$ are i.i.d. uniform in $[0, 1]$. The same is true when conditioning instead on $\widetilde{U}$, $I^U(\infty)$ and $A \subset I^U(\infty)$, since $I^U(\infty) = I^{\widetilde{U}}(\infty)$ is $\widetilde{U}$-measurable on the event $A \subset I^U(\infty)$. By averaging with respect to $\widetilde{U}$ and using $Z_A^U = Z_A^{\widetilde{U}}$, we conclude that, conditional on $I^U(\infty)$, on the event $A \subset I^U(\infty)$, the random variables $\langle U(e)/Z_\infty(e); e \in A \rangle$ are i.i.d. uniform in $[0, 1]$. The lemma follows.    $\square$

We also need the following fact:

LEMMA 3.19.    *If a graph $H$ of bounded degree does not contain a simple bi-infinite path, then $p_c(H) = 1$.*

PROOF.    By repeated applications of Menger's theorem, we see that if $x$ is a vertex in $H$, then there are infinitely many vertices $v$ such that $x$ is in a finite component of $H \setminus \{v\}$. Since $H$ has bounded degree, it follows that $p_c(H) = 1$.    $\square$

PROOF OF THEOREM 3.16.    Let $\omega_p$ be the set of edges $e$ satisfying $U(e) < p Z_\infty(e)$. Lemma 3.17 implies that, given $I(\infty)$, $\omega_p$ has the law of Bernoulli($p$) percolation on $I(\infty)$. Suppose that $\mathbf{P}[p_c(I(\infty)) < 1] > 0$ and fix $p < 1$ so that $\mathbf{P}[p_c(\omega_p) < 1] = \mathbf{P}[p_c(I(\infty)) < p] > 0$. On the event $p_c(\omega_p) < 1$, Lemma 3.19 implies that a.s., $\omega_p$ contains a simple bi-infinite path $\mathcal{P}$. Let $\alpha := \sup_{e \in \mathcal{P}} U(e)$. Since $\mathcal{P} \subset \omega_p$, we have, for $e \in \mathcal{P}$,

$$U(e) < p Z_\infty(e) \le p \sup_{\mathcal{P}} U = p\alpha,$$

which yields that $\alpha \le p\alpha$. Thus, $\alpha = 0$, which is clearly impossible. This contradiction shows that $p_c(I(\infty)) = 1$.    $\square$

Theorem 3.16 was stated for bounded degree graphs. If $G$ is locally finite but has unbounded degree, a similar argument still shows the following:

THEOREM 3.20.    *Let $G = (\mathsf{V}, \mathsf{E})$ be an infinite graph. Then the WMSF $\mathfrak{F}_w$ satisfies $p_c(\mathfrak{F}_w) = 1$ a.s., and moreover, $\bigcup_{v \in \mathsf{V}} I(v)$ has $p_c = 1$.*

PROOF. To see this, replace $Z_\infty(e)$ used in the above proof by $\widetilde{Z_\infty}(e)$, defined as $Z_\infty(e)$, except that the infimum ranges over all bi-infinite paths $\mathcal{P} \supset \{e\}$ that are edge-simple (no edge is repeated). □

COROLLARY 3.21. *Let $G$ be a unimodular transitive locally finite connected graph. Then $p_c(G) < p_u(G)$ iff $p_c(\mathfrak{F}_f) < 1$ a.s.*

PROOF. By Proposition 3.6 and Theorem 3.9, we have $p_c(G) < p_u(G)$ iff WMSF ≠ FMSF. Now the conclusion follows from Theorem 3.16 and Proposition 3.5(e). □

A dual argument shows that the FMSF is almost connected in the following sense:

THEOREM 3.22. *Let $G$ be any locally finite connected graph and $\varepsilon \in (0, 1)$. Let $\mathfrak{F}_f$ be a configuration of the FMSF and $\omega$ be an independent copy of $G[\varepsilon]$. Then $\mathfrak{F}_f \cup \omega$ is connected a.s.*

For this, we use a lemma dual to Lemma 3.17; it will provide a coupling of $\mathfrak{F}_f$ and $\mathfrak{F}_f \cup \omega$.

LEMMA 3.23. *Let $G = (\mathsf{V}, \mathsf{E})$ be a locally finite infinite graph and $\langle U(e); e \in \mathsf{E} \rangle$ be i.i.d. uniform $[0, 1]$ random variables. Conditioned on $\mathfrak{F}_f$, the random variables*

$$\frac{1 - U(e)}{1 - Z_f(e)} \qquad (e \notin \mathfrak{F}_f)$$

*are i.i.d. uniform $[0, 1]$ random variables.*

PROOF. Let $A \subset \mathsf{E}$ be a finite set such that $\mathbf{P}[A \cap \mathfrak{F}_f = \varnothing] > 0$. Let $\widetilde{U}(e) := 1$ for $e \in A$ and $\widetilde{U}(e) := U(e)$ for $e \notin A$, and let $Z_A^U$ denote the restriction of $Z_f^U$ to $A$. Consider any cut $\mathcal{C}$. If $A \cap \mathfrak{F}_f = \varnothing$ and $e \in A \cap \mathcal{C}$, then $U(e) > Z_f^U(e) \geq \inf\{U(e'); e' \in \mathcal{C} \setminus \{e\}\}$, by (3.1). Hence, if $A \cap \mathfrak{F}_f = \varnothing$, then for every cut $\mathcal{C}$,

$$\inf_{\mathcal{C}} U = \inf_{\mathcal{C} \setminus A} U = \inf_{\mathcal{C} \setminus A} \widetilde{U} = \inf_{\mathcal{C}} \widetilde{U}$$

and, therefore, (still assuming that $A \cap \mathfrak{F}_f = \varnothing$) $Z_A^U = Z_A^{\widetilde{U}}$ and $\mathfrak{F}_f^U = \mathfrak{F}_f^{\widetilde{U}}$, because $f \in \mathfrak{F}_f^U$ if and only if there is some cut $\mathcal{C} \ni f$ with $U(f) = \inf_{\mathcal{C}} U$, as we have seen in the proof of Lemma 3.2. Hence, $A \cap \mathfrak{F}_f^U = \varnothing$ implies $U > Z_A^{\widetilde{U}}$ on $A$. In fact, $A \cap \mathfrak{F}_f^U = \varnothing$ is equivalent to $U > Z_A^{\widetilde{U}}$ on $A$, because $Z_A^{\widetilde{U}} \geq Z_A^U$. Thus (by the observation in Remark 3.18), conditioned on $\widetilde{U}$ and $A \cap \mathfrak{F}_{\mathrm{f}}^{U} = \varnothing$, the random variables $\langle (1 - U(e))/(1 - Z_A^{\widetilde{U}(e)}); e \in A \rangle$ are i.i.d. uniform in $[0, 1]$. The same is true when conditioning instead on $\widetilde{U}$, $\mathfrak{F}_{\mathrm{f}}^{U}$ and $A \cap \mathfrak{F}_{\mathrm{f}}^{U} = \varnothing$, since $\mathfrak{F}_{\mathrm{f}}^{U} = \mathfrak{F}_{\mathrm{f}}^{\widetilde{U}}$ is $\widetilde{U}$-measurable on the event $A \cap \mathfrak{F}_{\mathrm{f}}^{U} = \varnothing$. By averaging with respect to $\widetilde{U}$ and using $Z_A^U = Z_A^{\widetilde{U}}$, we conclude that, conditional on $\mathfrak{F}_{\mathrm{f}}^{U}$, on the event $A \cap \mathfrak{F}_{\mathrm{f}}^{U} = \varnothing$, $\langle (1 - U(e))/(1 - Z_{\mathrm{f}}(e)); e \in A \rangle$ are i.i.d. uniform $[0, 1]$ variables. The lemma follows. $\quad\square$

PROOF OF THEOREM 3.22. By invoking Lemma 3.23, we see that $\mathfrak{F}_{\mathrm{f}} \cup \omega$ has the same law as

$$(3.4) \qquad\qquad \xi := \{e; 1 - U(e) \geq (1 - \varepsilon)[1 - Z_{\mathrm{f}}(e)]\}.$$

Thus, it suffices to show that $\xi$ is connected a.s. Consider any nonempty cut $\mathcal{C}$ in $G$, and let $\alpha := \inf_{e \in \mathcal{C}} U(e)$. Then $1 - \alpha = \sup_{\mathcal{C}} (1 - U)$, so we may choose $e \in \mathcal{C}$ to satisfy $1 - U(e) > (1 - \varepsilon)(1 - \alpha)$. By (3.1), $Z_{\mathrm{f}}(e) \geq \inf_{\mathcal{C} \setminus \{e\}} U \geq \alpha$ and, therefore, $e \in \xi$. Since $\xi$ intersects every nonempty cut, it is connected. $\square$

Recall that both WSF and WMSF have expected degree 2 in every transitive unimodular graph, by Theorem 3.12. We may combine Theorem 3.22 with recent work of Gaboriau [17] to deduce an inequality between the expected degrees of FSF and FMSF:

COROLLARY 3.24. Let $G = (\mathsf{V}, \mathsf{E})$ be a transitive unimodular connected infinite graph of degree $d$ and let $o \in \mathsf{V}$. Then

$$\mathbf{E}[\deg_o \mathrm{FSF}] \leq \mathbf{E}[\deg_o \mathrm{FMSF}] \leq 2 + d \int_{p_c}^{p_u} \theta(p)^2 \, dp,$$

where $\deg_v H$ denotes the degree of a vertex $v$ in a graph $H \ni v$.

The first inequality strengthens an observation of Lyons [29] (see the discussion following Conjecture 3.8 there), as well as its extension, Corollary 4.5, by Gaboriau [17]. The second inequality strengthens the inequality of Corollary 4.5 by Gaboriau [17].

PROOF OF COROLLARY 3.24. If $G$ is amenable, then FSF = WMSF and FMSF = WMSF by Corollary 6.3 of [10] and Proposition 3.5(c), so all have expected degree 2. Thus, the conclusion is trivial.

Assume now that $G$ is not amenable. We need some concepts defined or reviewed in [17]. First, there is a number $\beta_1(G) \geq 0$ called the first $\ell^2$-Betti number of $G$. Second, Theorems 6.4 and 7.8 of [10], together with Definition 2.9 of Gaboriau [17], show that $\mathbf{E}[\deg_o \mathrm{FSF}] = \mathbf{E}[\deg_o \mathrm{WSF}] + 2\beta_1(G) = 2 + 2\beta_1(G)$, an identity first observed by Lyons [30] for Cayley graphs. Third, let $\mathcal{O}_{\mathrm{HD}}$ be the set of graphs with no nonconstant harmonic Dirichlet functions. (A harmonic Dirichlet function on $G$ is a real-valued function $f$ on the vertex set of $G$ with $f(x) = \sum_{y \sim x} f(y)$ for all $x$ (harmonic) and $\sum_{x \sim y} [f(x) - f(y)]^2 < \infty$ (Dirichlet).) By Theorem 7.3 of [10], $\mathcal{O}_{\mathrm{HD}}$ consists precisely of the graphs on which the wired and free uniform spanning forests agree. By Theorem 4.2 of [17], if $\mu$ is an Aut(G)-invariant coupling of processes $\omega_1, \omega_2 \in 2^{\mathsf{E}}$ such that all clusters of $\omega_1$ are in $\mathcal{O}_{\mathrm{HD}}$ and $(\mathsf{V}, \omega_2)$ is connected, with both statements holding $\mu$-a.s., then

$$(3.5) \qquad 2\beta_1(G) \leq \sum_e \mu[e \in \omega_2 \setminus \omega_1],$$

where the summation is over edges incident with $o$. Let $\varepsilon > 0$ and let $\mu$ be the law of $(\mathfrak{F}_{\mathrm{w}}, \xi)$, where $\xi$ is defined in (3.4). Since $G$ is not amenable, we have $\theta(p_{\mathrm{c}}) = 0$, and so a.s. every cluster of $\mathfrak{F}_{\mathrm{w}}$ is a tree with one end, by Theorem 3.12. This implies that a.s. every cluster of $\mathfrak{F}_{\mathrm{w}}$ is in $\mathcal{O}_{\mathrm{HD}}$ (since the wired uniform spanning forest of a tree with one end is necessarily the whole tree). Combining this with Theorem 3.22, which shows that $\xi$ is a.s. connected with expected degree at most $\mathbf{E}[\deg_o \mathrm{FMSF}] + d\varepsilon$, we may apply (3.5) to this choice of $\mu$, obtaining $2\beta_1(G) \leq \mathbf{E}[\deg_o \mathrm{FMSF}] - \mathbf{E}[\deg_o \mathrm{WMSF}] + d\varepsilon = \mathbf{E}[\deg_o \mathrm{FMSF}] - 2 + d\varepsilon$, by Theorem 3.12. Since this holds for all positive $\varepsilon$, we have proved that $\mathbf{E}[\deg_o \mathrm{FSF}] \leq \mathbf{E}[\deg_o \mathrm{FMSF}]$.

To prove the final inequality, recall the event $A(e)$ used in the proof of Proposition 3.6; there, if $e$ is an edge with endpoints $x$ and $y$, then $A(e)$ was the event that $x$ and $y$ are in distinct infinite components of $(G \setminus \{e\})[U(e)]$. We saw there that $\mathbf{P}[A(e)] = \mathbf{P}[e \in \mathfrak{F}_{\mathrm{f}} \setminus \mathfrak{F}_{\mathrm{w}}]$. Let $A_x$ and $A_y$ be the events that $x$ and $y$ belong to infinite clusters, respectively. Since $\theta(p, G \setminus \{e\}, x) \leq \theta(p)$, the BK inequality of van den Berg and Kesten [42] gives $\mathbf{P}[A(e)] = \int_{p_{\mathrm{c}}}^{p_{\mathrm{u}}} \mathbf{P}[A(e) | U(e) = p] \, dp \leq \int_{p_{\mathrm{c}}}^{p_{\mathrm{u}}} \theta(p)^2 \, dp$. Sum this over all edges incident to $o$ to obtain the desired inequality. $\square$

## 4. Corollaries for planar graphs.

A *plane* graph is a graph $G$ embedded in the plane in such a way that no two edges cross each other. A *face* of a plane graph $G$ is a component of $\mathbb{R}^2 \setminus G$. A plane graph is *proper* if every bounded set in the plane contains only finitely many edges and vertices.

Suppose that $G$ is a proper plane graph. We define the *dual graph* $G^\dagger$ as follows. In each face $f$ of $G$, we place a single vertex $f^\dagger$ of $G^\dagger$. For every edge $e$ in $G$, we place an edge $e^\dagger$ in $G^\dagger$ connecting $f_1{}^\dagger$ and $f_2{}^\dagger$, where $f_1$ and $f_2$ are the two faces on either side of $e$. (It may happen that $f_1 = f_2$; then $e^\dagger$ is a loop.) This is done in such a way that $G^\dagger$ is a plane graph and $e \cap e^\dagger$ is a single point for every edge $e$ of $G$. Note that $G^\dagger$ is locally finite iff the boundary of every face of $G$ has finitely many vertices.

When $\Gamma$ is a subset of the edges of a plane graph $G$, define

$$\Gamma^* := \{e^\dagger ; e \notin \Gamma\}.$$

THEOREM 4.1 (FMSF is dual to WMSF). *Let $G$ and $G^\dagger$ be proper locally finite dual plane graphs. For any injection $U : \mathsf{E} \to \mathbb{R}$, let $U(e^\dagger) := 1 - U(e)$. We have*

$$(\mathfrak{F}_{\mathrm{f}}(U, G))^* = \{e^\dagger; U^\dagger(e^\dagger) < Z_{\mathrm{w}}^{U^\dagger}(e^\dagger)\},$$

*whence $(\mathfrak{F}_{\mathrm{f}}(U, G))^* = \mathfrak{F}_{\mathrm{w}}(U^\dagger, G^\dagger)$ if $U^\dagger(e^\dagger) \neq Z_{\mathrm{w}}^{U^\dagger}(e^\dagger)$ for all $e^\dagger \in \mathsf{E}^\dagger$.*

PROOF. The Jordan curve theorem implies that a set $\mathcal{P} \subset \mathsf{E} \setminus \{e\}$ is a simple path between the endpoints of $e$ iff the set $\mathcal{C} := \{f^\dagger; f \in \mathcal{P}\} \cup \{e^\dagger\}$ is a finite cut. Thus, the result is an immediate consequence of Lemma 3.2. □

The following easy result is proved in the same way that Proposition 12.5 of [10] is proved:

PROPOSITION 4.2 (Topology from duality). *Let $G$ be a proper plane graph with $G^\dagger$ locally finite. If each tree of the WMSF of $G$ has only one end a.s., then the FMSF of $G^\dagger$ has only one tree a.s. If, in addition, the WMSF of $G$ has infinitely many trees a.s., then the tree of the FMSF of $G^\dagger$ has infinitely many ends a.s.*

COROLLARY 4.3 ([6]). *The minimal spanning forest of $\mathbb{Z}^2$ is a.s. a tree with one end.*

PROOF. The hypothesis $\theta(p_{\mathrm{c}}) = 0$ of Theorem 3.12 applies by Harris [26] and Kesten [27]. Therefore, each tree in the WMSF has one end. By Proposition 4.2, this means that the FMSF has one tree. On the other hand, the wired and free measures are the same by amenability. □

The same reasoning shows the following:

PROPOSITION 4.4. *Let $G$ be a connected nonamenable proper plane graph with one end, such that there is a group of homeomorphisms of the plane acting quasi-transitively on the vertices of $G$. Then the FMSF on $G$ is a.s. a tree.*

The nonamenability assumption can be replaced by the assumption that the planar dual of $G$ satisfies $\theta(p_{\mathrm{c}}) = 0$. The latter assumption is known to hold in many amenable cases (see [28]).

PROOF OF PROPOSITION 4.4. Let $G$ be such an embedded graph. It is shown by Lyons with Peres [31] that $G$ and $G^\dagger$ have unimodular automorphism groups; the transitive case appeared in [13]. Thus, we may apply the main result of [8] to $G^\dagger$ to see that $\theta(p_c, G^\dagger) = 0$. Thus, Theorem 3.12 and Proposition 4.2 yield the desired conclusion. $\square$

We may also use similar reasoning to give another proof of a theorem of Benjamini and Schramm [13] (as extended from the transitive case by Lyons with Peres [31]):

COROLLARY 4.5. *If $G$ is a proper planar connected nonamenable graph with one end, such that there is a group of homeomorphisms of the plane acting quasi-transitively on the vertices of $G$, then $p_c(G) < p_u(G)$. In addition, there is a unique infinite cluster in Bernoulli($p_u(G)$) bond percolation.*

PROOF. By Theorem 3.9 and Proposition 3.6, it suffices to show that WMSF $\neq$ FMSF on $G$. Indeed, if they were the same, then they would also be the same on the dual graph, so that each would be one tree with one end, as in the proof of Proposition 4.4. But this is impossible by Theorem 5.3 of [9]. (Actually, that theorem was stated only in the transitive case, but the proof extends to quasi-transitive graphs.)

Furthermore, by Proposition 4.4, the FMSF is a tree on $G$, whence, by Proposition 3.13, there is a unique infinite cluster in Bernoulli($p_u$) percolation on $G$. $\square$

## 5. Correlations for the minimal spanning tree.

In view of the numerous similarities with uniform spanning trees, one might expect that the minimal spanning tree measure has negative associations. However, this is far from true. Indeed, the presence of even two edges can be positively correlated. To see this, we first present the following formula for computing probabilities of spanning trees. Let MST denote the minimal spanning tree measure on a finite connected graph.

PROPOSITION 5.1. *Let $G$ be a finite connected graph. Given a set $F$ of edges, let $N(F)$ be the number of edges of $G$ that do not become loops when each edge in $F$ is contracted. Note that $N(\varnothing)$ is the number of edges of $G$ that are not loops. Let $N'(e_1, \ldots, e_k) := \prod_{j=0}^{k-1} N(\{e_1, \ldots, e_j\})$. Let $T = \{e_1, \ldots, e_n\}$ be a spanning tree of $G$. Then*

$$\mathrm{MST}(T) = \sum_{\sigma \in S_n} N'(e_{\sigma(1)}, \ldots, e_{\sigma(n)})^{-1},$$

*where $S_n$ is the group of permutations of $\{1, 2, \ldots, n\}$.*

PROOF. To make the dependence on $G$ explicit, we write $N(F) = N(G; F)$. Note that $N(G/F; \varnothing) = N(G; F)$, where $G/F$ is the graph $G$ with each edge in $F$ contracted. Given an edge $e$ that is not a loop, the chance that $e$ is the least edge in the minimal spanning tree of $G$ equals $N(G; \varnothing)^{-1}$. Furthermore, given that this is the case, the ordering on the nonloops of the edge set of $G/\{e\}$ is uniform. Thus, if $f$ is not a loop in $G/\{e\}$, then the chance that $f$ is the next least edge in the minimal spanning tree of $G$, given that $e$ is the least edge in the minimal spanning tree of $G$, equals $N(G/\{e\}; \varnothing)^{-1} = N(G; \{e\})^{-1}$. Thus, we may easily condition, contract and repeat.

Thus, the probability that the minimal spanning tree is $T$ *and* that $e_{\sigma(1)} < \cdots < e_{\sigma(n)}$ is equal to

$$\prod_{j=0}^{n-1} N(G/\{e_{\sigma(1)}, e_{\sigma(2)}, \ldots, e_{\sigma(j)}\}; \varnothing)^{-1} = N'(e_{\sigma(1)}, \ldots, e_{\sigma(n)})^{-1}.$$

Summing this over all possible induced orderings of $T$ gives $\mathrm{MST}(T)$.  $\square$

An example of a graph where the inclusion of two specific edges in the MST is positively correlated is provided by the complete graph $K_4$, with two of its edges that do not share endpoints replaced by three edges (each) in parallel. If $e_1$ and $e_2$ are two of these parallel edges not sharing endpoints with each other, then $\mathrm{MST}[e_1, e_2 \in T] > \mathrm{MST}[e_1 \in T]\mathrm{MST}[e_2 \in T]$; the left-hand side divided by the right-hand side turns out to be $109872/109561$. To aid the reader who wishes to check the calculations, we present the following outline, where the probabilities that appear below were calculated using Proposition 5.1. Let $e_3$, $e_4$, $e_5$ and $e_6$ be the four edges that are not replaced by parallel ones, with $e_3$ and $e_4$ not incident with each other and $e_1, e_3, e_5$ sharing a common vertex. Let us use the following shorthand for a spanning tree: $[ijk]$ will denote the tree formed by the edges $e_i, e_j, e_k$. Then there are four types of spanning trees, where "isomorphic" below refers to the existence of an automorphism of $G$ that sends one spanning tree to another:

- 12 trees isomorphic to $[134]$, each with probability $163/12600$;
- 36 trees isomorphic to $[123]$, each with probability $109/6300$;
- 4 trees isomorphic to $[345]$, each with probability $7/600$;
- 12 trees isomorphic to $[135]$, each with probability $23/1575$.

Thus,  $\mathrm{MST}[e_1 \in T] = \mathrm{MST}[e_2 \in T] = 331/1260$  and  $\mathrm{MST}[e_1, e_2 \in T] = 109/1575$.

**6. Concluding remarks and questions.** We give some examples to illustrate how the behavior of minimal spanning forests may differ from that of uniform spanning forests, or how we can answer certain questions for minimal spanning forests that are still open for uniform spanning forests. We then give some open questions and conjectures.

EXAMPLE 6.1. We describe a connected planar graph where the FMSF equals the WMSF and has two components a.s. If we drop the planarity requirement, then such an example is known for the uniform spanning forests (Remark 9.8 in [10]). That example can be modified so as to be planar, but with unbounded degrees; a bounded-degree example is impossible by Corollary 12.9 of [10]. Let $G_-$ be the lower half plane of $\mathbb{Z}^2$, that is, the subgraph of $\mathbb{Z}^2$ spanned by vertices $(x, y) \in \mathbb{Z}^2$ with $y \leq 0$. For $n > 0$, let $\tau(n)$ be the probability that invasion percolation on $G_-$ from the origin reaches $(n, 0)$. By Barsky, Grimmett and Newman [7], in $G_-$ we have $\theta(p_c) = 0$, and for every $p > p_c$, there is a unique infinite $G_-[p]$ cluster. Moreover, for $p > p_c$, a.s. every invasion tree meets the infinite $G_-[p]$ cluster; see [24], Proposition 4.3.1. We claim that $\lim_{n \to \infty} \tau(n) = 0$. For any fixed $p > p_c$, the event that $(n, 0)$ is in the invasion tree of the origin is contained in the union of two events: the event that $(n, 0)$ is in an infinite $G_-[p]$ cluster, and the event that the invasion tree of the origin has not met an infinite $G_-[p]$ cluster in its first $n$ invasion steps. The probability of the latter event tends to zero as $n \to \infty$ and the probability of the former event tends to 0 as $p \downarrow p_c$, because $\theta(p)$ is right continuous at $p_c$. Consequently, $\lim_{n \to \infty} \tau(n) = 0$, as claimed. Thus, we may choose a sequence $\langle n_k \rangle$ such that $\sum_k \tau(n_k) < \infty$. Let $G$ be the subgraph of $\mathbb{Z}^2$ obtained by erasing all edges connecting $(n, 0)$ to $(n, 1)$ for $n \notin \{n_k\}$ (in particular, all these edges for $n \leq 0$ are erased). Since $p_c$ of the lower half plane is the same as $p_c$ of the whole plane, as shown by Grimmett and Marstrand [19], we have $p_c(G) = p_c(\mathbb{Z}^2) = 1/2$. Since critical percolation in $\mathbb{Z}^2$ has infinitely many disjoint open cycles that surround the origin (see [18]), there are infinitely many disjoint paths in the lower half plane that are open at level $1/2$ and connect the negative $x$-axis to the positive $x$-axis. A similar situation obtains for the upper half plane. For any vertex $v$ in $G$, the invasion tree $T_v$ must intersect, and thus fill, infinitely many of the paths in one of these two families. It follows that the WMSF in $G$ has at most two components.

Given $\varepsilon > 0$, choose $m > 0$ such that $\sum_k \tau(m + n_k) < \varepsilon$. Then the vertices $(-m, 0)$ and $(-m, 1)$ have disjoint invasion trees with probability at least $1 - 2\varepsilon$. This shows that the WMSF in $G$ a.s. has two components, $\Gamma_1$ and $\Gamma_2$. We claim that these are also the components of the FMSF. Indeed, let $e$ be an edge with endpoints in $\Gamma_1$ and $\Gamma_2$. Clearly, $U(e) > 1/2$. Let $p \in (1/2, U(e))$ be rational. By Barsky, Grimmett and Newman [7], there is a unique Bernoulli($p$) percolation infinite cluster in each half plane, and the vertices $(n_k, 0)$ and $(n_k, 1)$ are in these clusters for infinitely many $k$. Therefore, there is a unique Bernoulli($p$) percolation cluster in $G$. The invasion trees from the endpoints of $e$ intersect this cluster since $\theta(1/2, G) = 0$, and the labels on these trees are all less than $U(e)$ (otherwise, $e$ would be invaded). Therefore, $e \notin$ FMSF. Thus, the FMSF equals the WMSF on $G$.

The number of trees in the WSF is an a.s. constant on every graph. This is not true for the WMSF:

EXAMPLE 6.2.   Next, we provide a planar graph where the number of trees in the WMSF is not an a.s. constant. Let $H_+$ be the upper half of $\mathbb{Z}^2$ (above the $x$-axis). Let $p_k := 1/2 + 2^{-k-1}$. Choose $\{m_k\}$ so that (i) with probability greater than $1 - k^{-2}$, an infinite cluster of Bernoulli($p_k$) percolation in $H_+$ comes within distance $m_k$ of the origin, and (ii) with probability greater than $1 - k^{-2}$, Bernoulli($1/2$) percolation in $\mathbb{Z}^2$ has an open cycle in the annulus $B(0, m_{k+1}) \setminus B(0, m_k)$. Define a planar graph $G_+$ by adding to every edge in $H_+$ a parallel path joining its endpoints; if the distance from the edge to the $x$-axis is between $m_k$ and $m_{k+1}$, then the path has length $k$. It is easy to check that $p_c(G_+) = 1/2$, that there is a unique infinite cluster for Bernoulli($1/2$) percolation on $G_+$ a.s. and that every invasion tree in $G_+$ intersects that cluster a.s. Let $G_-$ be a copy of $G_+$, built on the lower half plane. Let $v_+$ and $v_-$ be two fixed vertices on the boundaries of $G_+$ and $G_-$, respectively. Add a single edge $e := [v_+, v_-]$, thus defining a planar graph $G$ with $p_c(G) = 1/2$. Consider the event $A$ that each of $v_+$, $v_-$ contained in the unique critical infinite cluster in $G_+$, respectively $G_-$. If $U(e) > 1/2$ and $A$ occurs, then the WMSF has two trees. If $U(e) < 1/2$ and $A$ occurs, then the WMSF has one tree.

It was asked in [10] whether the number of trees in the FSF is always an a.s. constant. This is still not known, but we can show it is not true for the FMSF:

EXAMPLE 6.3.   The previous example can be modified to obtain a planar graph $\Gamma$ where the number of trees in the FMSF is not an a.s. constant. Define $\Gamma$ by joining $G_+$ and $G_-$ at boundary vertices $(j, \pm 1)$ with $\lceil 2^j/j^2 \rceil$ disjoint paths of length $j$, for each $j \geq 1$. Consider the infinite clusters $\Delta_+$ and $\Delta_-$ for Bernoulli($1/2$) percolation in $G_+$, respectively $G_-$. For any $p > 1/2$, infinitely many of the paths we added will have all labels less than $p$ by Borel–Cantelli; therefore, there is a.s. a unique infinite cluster for Bernoulli($p$) percolation in $\Gamma$.

If there is exactly one path $\mathcal{P}$ in $\Gamma$ connecting $\Delta_+$ and $\Delta_-$ with $U(e) < 1/2$ for all $e \in \mathcal{P}$, then the FMSF in $\Gamma$ is connected, because the same holds for the WMSF, by the analysis of the preceding example; if there is no such path $\mathcal{P}$, then the FMSF in $\Gamma$ has two components. Both of these possibilities occur with positive probability, by Borel–Cantelli.

EXAMPLE 6.4.   Although for transitive unimodular graphs, the inequality WSF $\neq$ FSF implies the inequality WSF $\neq$ FSF (see Corollary 3.24 and recall that the expected degrees of the vertices in the WSF and in the WMSF are 2), the same is not true for general graphs. We give an example to illustrate this phenomenon: Define $\tau_p(r)$ to be the probability that Bernoulli($p$) percolation on $\mathbb{Z}^3$ connects the origin to any point at distance at least $r$ from the origin. Given a positive integer $m$, write $G_m$ for the graph obtained from $\mathbb{Z}^2$ by replacing each edge by $m$ parallel edges. Choose and fix $m$ so that $p_c(G_m) < p_c(\mathbb{Z}^3)$. Let $G_{m,n}$ be the subgraph of $G_m$ induced by an $n$-by-$n$ square. Write $C_n$ for the effective conductance between one corner of $G_{m,n}$ and the opposite corner. (All edges are given unit conductance.) Because $\mathbb{Z}^2$, and hence $G_m$, is recurrent, $C_n \to 0$ as $n \to \infty$. Choose a sequence $\langle n_k \rangle$ such that

$$(6.1) \qquad \sum_k C_{n_k} < \infty.$$

Also, choose an increasing sequence $\langle R_k \rangle$ such that

$$(6.2) \qquad \sum_k \tau_{p_c(\mathbb{Z}^3)-1/k}(R_k) < \infty,$$

and fix a sequence of vertices $\langle x_k \rangle$ in $\mathbb{Z}^3$ such that $\|x_k - x_j\| \geq R_k$ for all pairs $j < k$.

Finally, let $G$ be the graph obtained by starting with two copies of $\mathbb{Z}^3$ and identifying $x_k$ in one copy with one corner of $G_{m,n_k}$ and identifying $x_k$ in the other copy with the opposite corner of $G_{m,n_k}$. We claim that $\mathrm{WSF}(G) \neq \mathrm{FSF}(G)$, while $\mathrm{WMSF}(G) = \mathrm{FMSF}(G)$.

To prove this, we first show that $p_c(G) = p_c(\mathbb{Z}^3)$. Indeed, since $G$ contains a copy of $\mathbb{Z}^3$, we have trivially that $p_c(G) \leq p_c(\mathbb{Z}^3)$. On the other hand, given $p < p_c(\mathbb{Z}^3)$, condition (6.2) and the Borel–Cantelli lemma ensure that in each copy of $\mathbb{Z}^3$, a.s. all but finitely many $p$-clusters contain at most one of the points $x_k$, and the other $p$-clusters contain finitely many of the points $x_k$. Therefore, all $p$-clusters of $G$ are finite a.s.

Next, note that for $p > p_c(G)$, we have a.s. that infinitely many $k$ are such that each copy of $x_k$ lies in the infinite $p$-cluster of its copy of $\mathbb{Z}^3$ and also, since $p_c(G) > p_c(G_m)$, both copies of $x_k$ are connected to each other in $G_{m,n_k}$. It follows that there is a.s. a unique infinite $p$-cluster in $G$, whence, by Proposition 3.6, we obtain that $\mathrm{WMSF}(G) = \mathrm{FMSF}(G)$.

It remains to show that $\mathrm{WSF}(G) \neq \mathrm{FSF}(G)$. For this, it suffices (and indeed, is equivalent) to show that there is a nonconstant harmonic Dirichlet function on $G$; see [10], Theorem 7.3 for this criterion. We now define such a function. Let $f(x)$ be the probability that when simple random walk starts from $x$ in $G$, it eventually stays in the first copy of $\mathbb{Z}^3$. To show that $f$ is not constant, take a starting point that is very far from $x_1, x_2, \ldots, x_K$. The expected number of visits the walk makes to $x_k$ for $1 \leq k \leq K$ is then very small and the expected number of visits to $x_k$ for $k > K$ is bounded. On each visit to $x_k$, the chance of crossing $G_{m,n_k}$ to the other copy of $x_k$ before leaving $G_{m,n_k}$ is at most $aC_{n_k}$ for some constant $a$ (see, e.g., [15] or [31]). From (6.1), it follows that the expected number of crossings from one copy of $\mathbb{Z}^3$ to the other is very small, and hence the probability of making any crossing is also very small. That is, $f(x)$ tends to 1 as $x \to \infty$ in the first copy of $\mathbb{Z}^3$, while $f(x)$ tends to 0 as $x \to \infty$ in the second copy of $\mathbb{Z}^3$. The fact that $f$ is harmonic is obvious from the Markov property of simple random walk. Finally, to see that $f$ is Dirichlet, observe that $f = \lim_n f_n$, where $f_n(x)$ is the probability that the random walk starting at $x$ will hit a vertex $v$ satisfying $\|v\| > n$ in the first copy of $\mathbb{Z}^3$ before hitting a vertex $v$ satisfying $\|v\| > n$ in the second copy. Let $g$ be the function that equals one in the first copy of $\mathbb{Z}^3$, zero in the second copy, and is harmonic at other vertices. The Dirichlet energy $\sum_{x \sim y}[g(x) - g(y)]^2$ is equal to $\sum_k C_{n_k}$, which is finite. Since on finite graphs harmonic functions with given boundary values minimize the Dirichlet energy, the Dirichlet energy of $f_n$ is no larger than that of $g$. The same therefore holds for $f = \lim_n f_n$.

There are many open questions related to minimal spanning forests:

QUESTION 6.5. Let $G$ be a transitive graph whose automorphism group is not unimodular. Does every tree of the WMSF on $G$ have one end a.s.?

QUESTION 6.6. Does every nonamenable quasi-transitive graph $G$ satisfy FMSF $\neq$ WMSF? In view of Proposition 3.6, a positive answer is equivalent to a well-known conjecture by Benjamini and Schramm [11] that any such graph $G$ satisfies $p_c(G) < p_u(G)$.

QUESTION 6.7. If $G$ is a unimodular transitive graph and WMSF $\neq$ FMSF, does FMSF-a.s. every tree have infinitely many ends?

After the first version of this paper was circulated, a proof of this conjecture was obtained by Timár [41].

QUESTION 6.8. For which $d$ is the minimal spanning forest of $\mathbb{Z}^d$ a.s. a tree? This question is due to Newman and Stein [35] who conjecture that the answer is $d < 8$ or $d \le 8$.

There is a related conjecture of Benjamini and Schramm [12]:

CONJECTURE 6.9. Let $G$ be a quasi-transitive nonamenable graph. Then FMSF is a single tree a.s. iff there is a unique infinite cluster in $G[p_u]$ a.s.

We can strengthen this conjecture to say that the number of trees in the FMSF equals the number of infinite clusters at $p_u$. An even stronger conjecture would be that in the natural coupling of Bernoulli percolation and the FMSF, each infinite cluster at $p_u$ intersects exactly one component of the FMSF and each component of the FMSF intersects exactly one infinite cluster at $p_u$. Recall that in Proposition 3.13, we proved the second part of this conjecture for those $G$ that satisfy $p_u(G) > p_c(G)$.

QUESTION 6.10. Must the number of trees in the FMSF and the WMSF in a quasi-transitive graph be either 1 or $\infty$ a.s.? This question for $\mathbb{Z}^d$ is due to Newman [34].

CONJECTURE 6.11. The components of the FMSF on a unimodular transitive graph are indistinguishable in the sense that for every automorphism-invariant property $\mathcal{A}$ of subgraphs, either a.s. all components satisfy $\mathcal{A}$ or a.s. all do not. The same holds for the WMSF.

This fails in the nonunimodular setting, as the example in [32] shows.

CONJECTURE 6.12. Let $T_o$ be the component of the identity $o$ in the WMSF on a Cayley graph, and let $\xi = \langle v_n; n \geq 0 \rangle$ be the unique ray from $o$ in $T_o$. The sequence of "bushes" $\langle b_n \rangle$ observed along $\xi$ converges in distribution. (Formally, $b_n$ is the connected component of $v_n$ in $T \setminus \{v_{n-1}, v_{n+1}\}$, multiplied on the left by $v_n^{-1}$.)

QUESTION 6.13. Given a finitely generated group $\Gamma$, does the expected degree of a vertex in the FMSF of a Cayley graph of $\Gamma$ depend on which Cayley graph is used? As discussed in the proof of Corollary 3.24, the analogous result is true for the FSF.

QUESTION 6.14. One may consider the minimal spanning tree on $\varepsilon \mathbb{Z}^2 \subset \mathbb{R}^2$ and let $\varepsilon \to 0$. It would be interesting to show that the limit exists in various senses. Aizenman, Burchard, Newman and Wilson [2] have shown that a subsequential limit exists.

**Acknowledgments.** We thank Itai Benjamini for useful discussions at the early stages of this work, and Gábor Pete for comments on the manuscript. Part of this work was done while the first two authors were visiting Microsoft Research.

R. LYONS
DEPARTMENT OF MATHEMATICS
INDIANA UNIVERSITY
BLOOMINGTON, INDIANA 47405-5701
USA
E-MAIL: rdlyons@indiana.edu
URL: http://mypage.iu.edu/~rdlyons/

Y. PERES
DEPARTMENTS OF STATISTICS AND MATHEMATICS
UNIVERSITY OF CALIFORNIA
BERKELEY, CALIFORNIA 94720-3860
USA
E-MAIL: peres@stat.berkeley.edu
URL: http://www.stat.berkeley.edu/~peres/

O. SCHRAMM
MICROSOFT RESEARCH
ONE MICROSOFT WAY
REDMOND, WASHINGTON 98052
USA
URL: http://research.microsoft.com/~schramm/